\documentclass[12pt,a4paper]{article}
\usepackage{amsthm,amssymb,amsmath,mathrsfs,graphicx}
\usepackage[pagebackref=false,colorlinks,linkcolor=black,citecolor=magenta]{hyperref}
\allowdisplaybreaks[1]

\newtheorem{theorem}{\bf Theorem}
\newtheorem{proposition}{\bf Proposition}
\newtheorem{lemma}{\bf Lemma}

\newtheorem{corollary}{\bf Corollary}
\newtheorem{remark}{\bf Remark}
\newtheorem{prob}{\bf Problem}
\newtheorem{defin}{\bf Definition}
\newtheorem{conj}{\bf Conjecture}

\renewcommand{\proof}{\noindent{\it\textbf{Proof.}}\ \ }

\newcommand{\diag}{{\rm diag}}

\newcommand{\Cay}{{\rm Cay}}

\newcommand{\eqd}{$\hfill \blacksquare$}

\newcommand{\Irr}{{\rm Irr}}

\date{}

\title{Integral Cayley graphs of symmetric groups on transpositions}

\author{Alireza Abdollahi $^{a}$ \footnote{E-mail:
a.abdollahi@math.ui.ac.ir (corresponding author)}, Majid Arezoomand $^{b}$ \footnote{Email: arezoomand@lar.ac.ir} and Mahdi Ebrahimi $^{c}$ \footnote{Email: m.ebrahimi.math@ipm.ir}\\
{\small\em$^a$ Department of Pure Mathematics, Faculty of Mathematics and Statistics,}\\ 
{\small\em University of Isfahan, Isfahan 81746-73441, Iran}\\
{\small\em$^b$ University of Larestan, Larestan 74317-16137, Iran}\\
{\small\em$^c$School of Mathematics, Institute for Research in Fundamental Sciences (IPM),}\\
{\small\em Tehran 19395-5746, Iran}
}

\begin{document}
\maketitle
\begin{abstract}
We study subsets $T$ consisting of some transpositions $(i,j)$ of the symmetric group $S_n$ on $\{1,\dots,n\}$ such that the Cayley graph $\Gamma_T:=Cay(S_n,T)$ is an integral graph, i.e., all eigenvalues of an adjacency matrix of $\Gamma_T$ are integers. 
Graph properties of  $\Gamma_T$ are determined in terms of ones of the graph $G_T$ whose vertex set is $\{1,\dots,n\}$ and $\{i,j\}$ is an edge if and only if $(i,j)\in T$. 
Here we prove that if $G_T$ is a tree then $\Gamma_T$ is integral if and only if $T$ is isomorphic to the star graph $K_{1,n-1}$, answering Problem 5 of [Electron. J. Comnin., 29(2) (2022) \# P2.9]. Problem 6 of the latter article asks to find necessary and sufficient conditions on $T$ for  integralness of $Cay(S_n,T)$ without any further assumption on $T$. We show that if $G_T$ is a graph which we  call it a ``generalized complete multipartite graph" then $Cay(S_n,T)$ is integral. We conjecture that $Cay(S_n,T)$ is integral only if $G_T$ is a generalized complete multipartitie graph. To support the latter conjecture we show its validity whenever $G_T$ is some classes of graphs including cycles and cubic graphs.        
\end{abstract}

\medskip \noindent
Keywords: Cayley graph; Symmetric group; Transposition; Spectra of graphs; Eigenvalues of graphs; Integral graphs.\\
Mathematics Subject Classification: 05C50; 05C25.

\section{Introduction}

Let $G$ be a group and $S$ be a subset of $G$ such that $S$ is inverse closed (i.e. $s\in S\Leftrightarrow s^{-1} \in S$) and the identity element $1$ of $G$ does not belong to $S$. Then the Cayley graph $Cay(G,S)$ is the graph whose vertex set is the set of elements of $G$ and two vertices $g_1,g_2\in G$ are adjacent if and only if $g_1 g_2^{-1} \in S$. It is well known that $Cay(G, S)$ is connected if and only if $G$ can be generated by $S$. 

Let $S_n$ be the symmetric group on $[n]:=\{1,\dots,n\}$.
For a set $T$ of transpositions in $S_n$, $G_T$ is the graph with vertex set $[n]$ such that $i, j \in [n]$ are adjacent if and only if $(i, j) \in T$. It is well known
that $Cay(S_n,T)$ is connected if and only if $G_T$ is connected. 
The (adjacency) spectrum of a finite simple graph $\Gamma$ is the multiset of eigenvalues of an adjacency matrix of $\Gamma$. A graph is called integral if and only if the spectrum consists of only integers.\\

Here we study the following problems proposed in \cite{LZ}.

\begin{prob}{\rm \cite[Problem 5]{LZ}} \label{p5}
 Let $T$ be a set of transpositions in $S_n$ such that $G_T$ is a tree, where $n \geq 2$.
Give a necessary and sufficient condition for $Cay(S_n,T)$ to be integral.
\end{prob}
\begin{prob}{\rm \cite[Problem 6]{LZ}} \label{p6} Let $T$ be a set of transpositions in $S_n$ such that $G_T$ is connected, where $n \geq 2$.
Give a necessary and sufficient condition for $Cay(S_n,T)$ to be integral.
\end{prob}

In Problem \ref{p5}, a sufficient condition on $T$ is known to be $G_T\cong K_{1,n-1}$ the star graph with $n$ vertices (see \cite{AV}, \cite[Theorem 5.3]{C}, \cite{CF}, \cite{KM}).  We complete the answer to Problem \ref{p5} by showing that the latter condition is also necessary (see part (b) of Corollary \ref{cut}, below).

About Problem \ref{p6}, we prove that if $G_T$ is a ``generalized complete multipartite graph" (see Definition \ref{def1}) then $Cay(S_n,T)$ is integral. We propose as a conjecture (see Conjecture \ref{conj}, below) that the converse of the latter holds valid. To support Conjecture \ref{conj}, we prove its validity for some classes of graphs including cycles and cubic graphs (see Propositions \ref{prop1}, \ref{prop2} and \ref{prop3}).   

\section{Preliminaries}

Throughout the paper we use the following notations:

$C_n$: the cycle graph with $n$ vertices.

$K_n$: the complete graph with $n$ vertices.

$K_{n_1,n_2,\dots,n_k}$: the complete multipartite graph with $k$ parts of sizes $n_1,n_2,\dots,n_k$.

For a given graph $\Gamma$:

$V(\Gamma)$:  the set of vertices of $\Gamma$.

$E(\Gamma)$:  the set of edges of $\Gamma$. 

$n\Gamma$: the graph which is disjoint union of $n$ copies of the graph $\Gamma$.

$\Gamma_1\cup\Gamma_2$: the disjoint union of graphs $\Gamma_1$ and $\Gamma_2$.

$\Gamma_1*\Gamma_2$: the join graph of disjoint graphs $\Gamma_1$ and $\Gamma_2$.

$\mathcal{L}(\Gamma)$: the Laplacian matrix of $\Gamma$ which is $D-A$, where $D$ is the diagonal matrix of the vertex-degrees and $A$ is the adjacency matrix of $\Gamma$ with the same row vertex-labeling as $D$.   

$I_\Gamma$: the set of all isolated vertices of a graph $\Gamma$.

$\Gamma^\#$ : the induced subgraph of $\Gamma$ on the set $V(\Gamma)\setminus I_\Gamma$.

$\Gamma+v$: the graph obtained by adding the isolated vertex $v$ to the graph $\Gamma$.

\subsection{Eigenvalues of Cayley graphs}

For the basic definitions and facts from representation theory of finite groups, not defined here, we refer the reader to \cite{Serre}. Let $G$ be a finite group and $\Irr(G)=\{\rho_1,\ldots,\rho_m\}$ be the set of all inequivalent irreducible complex representations of $G$ and $\Bbb CG$ be the group algebra of elements of the form $\sum_{g\in G}a_g g$, where $a_g\in\Bbb C$.
It is clear that $\Bbb CG$ is a vector space over $\Bbb C$ of dimension $|G|$. The \textit{left regular representation} $\rho_{\textrm{reg}}:G\rightarrow GL(\Bbb CG)$ is the right translation acting naturally on $\Bbb CG$:
\[\rho_{\textrm{reg}}(g)(\sum_{h\in G}a_h h)=\sum_{h\in G}a_h(gh).\]

It is well-know that $\rho_{\textrm{reg}}=d_1\rho_1\oplus\dots\oplus d_m\rho_m$, where $d_i$ is the degree of $\rho_i$.  Let $\Gamma=\Cay(G,S)$ be a Cayley (di)graph over $G$ with respect to a subset $S$ of $G$ and $A$ be the adjacency matrix of $\Gamma$. If we view $A$ as a linear map $A:\Bbb CG\rightarrow\Bbb CG$, then $A=\sum_{g\in G} \rho_{\textrm{reg}}(g)$ and the following result arises:

\begin{proposition}{\rm (\cite[Corollary 7]{AT} and \cite{DS})}\label{at}
Let $\Gamma={\rm Cay}(G,S)$ be an Cayley  graph over a
finite group $G$ with irreducible matrix representations $\rho_1,\ldots,\rho_m$.
Let $d_l$ be the degree of $\rho_l$. For each $l\in\{1,\ldots,m\}$, define a
$d_l\times d_l$  matrix $A_l:=\rho_l(S)$, where $\rho_l(S)=\sum_{s\in S}\rho_l(s)$. Let
$\chi_{A_l}(\lambda)$ and $\chi_{A}(\lambda)$ be the
characteristic polynomials of $A_l$ and $A$, respectively. Then
\begin{enumerate}
\item[(1)] there exists a basis $\mathcal{B}$ such that
$[A]_{\mathcal{B}}={\rm Diag}(A_1\otimes I_{d_1},\ldots,A_m\otimes I_{d_m})$, where $I_{d_l}$ is the $d_l \times d_l$ identity matrix.
\item[(2)]$\chi_A(\lambda)=\Pi_{l=1}^{m}\chi_{A_l}(\lambda)^{d_l}$.
\end{enumerate}
\end{proposition}

A partition is a weakly decreasing finite sequence of positive integers $\alpha=(\alpha_1, \dots, \alpha_l)$. We call $|\alpha|=\alpha_1+\dots+\alpha_l$ the \textit{size} of $\alpha$. The notation $\alpha \vdash n$ is used for a partition $\alpha$ of a positive integer $n$. It is well known that both the conjugacy classes of $S_n$ and the irreducible characters of $S_n$ are indexed by partitions $\alpha$ of $[n]$ (see \cite{GA}). We 
denote the irreducible character corresponding to the partition $\alpha\vdash n$ by $\chi^\alpha$ and its degree by $f_\alpha$.

\begin{lemma}\label{lemm1} Let $T$ be the set of all transpositions of $S_n$ or equivalently $G_T$ is the complete graph. Then  $\Cay(S_n,T)$ is an integral graph and its spectrum is 
$$\big\lbrace{q_\alpha:=\frac{n(n-1)}{2f_\alpha}\chi^\alpha((1,2)) \;|\; \alpha\vdash n\big\rbrace}.$$
\end{lemma}
\proof Since the set $T$ is a conjugacy class of $S_n$,  the first part follows from \cite[Lemma 5]{DS}.    Being integer of $q_\alpha$ follows from  the fact that the eigenvalues of a simple graph are algebraic integers and by \cite[Lemma 7]{DS} that $q_\alpha$ is a rational number. 
\eqd

\begin{proposition}\label{multipartite}
Let $T$ be a set of transpositions in $S_n$ such that $G_T$ is a disjoint union of complete multipartite graphs. Then $\Cay(S_n,T)$ is integral.
\end{proposition}
\proof Since by Lemma \ref{lemm1} for every partition $\alpha\vdash n$, $q_\alpha$ is an integer, it follows from \cite[Theorem 5.3 and Proposition 5.2]{C}, that $\Cay(S_n,T)$ is integral.
\eqd
\section{Main results}

Given a permutation $\pi\in S_n$, we denote by $P_\pi=(p_{ij})$ the $n\times n$ permutation matrix with entries $p_{ij}=1$ if $i^\pi=j$, and $0$ otherwise.

Recall that a graph $\Gamma$ is called Laplacian integral if the spectrum of $\mathcal{L}(\Gamma)$ consists of only integers. 

\begin{theorem}\label{main}
Let $T$ be a set of transpositions in $S_n$, $n\geq 2$. If $\Cay(S_n,T)$ is integral then $G_T$ is Laplacian integral.
\end{theorem}
\proof
Given $\tau=(i,j)\in T$, we have $P_\tau=I_n-e_{ii}-e_{jj}+e_{ij}+e_{ji}$, where $e_{rs}$ is the $n\times n$ matrix that its $rs$-th entry is $1$ and the remaining entries are $0$. On the other hand, $\sum_{(i,j)\in T}[e_{ii}+e_{jj}]$ and $\sum_{(i,j)\in T}[e_{ij}+e_{ji}]$ are the degree matrix and the adjacency matrix of $G_T$, respectively. Therefore,
$\sum_{\tau\in T}P_\tau=|T|I_n-\mathcal{L}(G_T)$. On the other hand, each eigenvalue of $\sum_{\tau\in T}P_\tau$ is an eigenvalue of $\Cay(S_n,T)$, by \cite[Proposition 2.1]{DF}. Therefore integrality of $\Cay(S_n,T)$ implies the integrality of
$\mathcal{L}(G_T)$ as desired.
\eqd
\begin{remark}\label{rmk} One can see that the map
\begin{eqnarray*}
\varphi&:&S_6\rightarrow GL_5(\Bbb C)\\
(1,2,3,4,5)&\mapsto & \begin{bmatrix}
 0& 0& 0& 0& 1\\
 1& 0& 0& 0& 0\\
 0& 0& 0& 1& 0\\
 0& 0& 1& 0& 0\\
 0& 1& 0& 0& 0
\end{bmatrix} \\
(1,2)&\mapsto & \begin{bmatrix}
0& 0& 1& 0& 0\\
0& 0& 0& 0& 1\\
1&0& 0& 0& 0\\
-1& -1& -1& -1& -1\\ 
0& 1& 0& 0& 0 
\end{bmatrix}
\end{eqnarray*}
is an irreducible representation of $S_6$. Let $T=\{(1,2),(2,3),(3,4),(4,5),(5,6),(6,1)\}$. Then the characteristic polynomial
of $\sum_{t\in T}\varphi(t)$ is $x(x+1)^2(x^2+2x-12)$ has non-integer roots, which means that $\Cay(S_n,T)$ is  non-integral Cayley graph by Proposition \ref{at}. Since $G_T\cong C_6$ is Laplacian integral, this means that the converse of 
Theorem \ref{main} does not hold in general.
\end{remark}

\begin{defin} \label{defin0} {\rm Suppose $G_1,\dots,G_k$, $k\geq 2$,  are  graphs with the same vertex set. A graph $G$ is called \textit{weakly join of $G_1,\dots,G_k$} if there exist subgraphs $\Gamma_i$ of $G_i$, $1\leq i\leq k$ such that 
\begin{itemize}
\item[(1)] for every $1\leq i\leq k$, $E(\Gamma_i)=E(G_i)$,
\item[(2)] for every $i\neq j$, $V(\Gamma_i)\cap V(\Gamma_j)=\varnothing$, and
\item[(3)] $G^\#=\Gamma_1*\dots*\Gamma_k$.
\end{itemize}
We use the notation $\mathcal{WJ}(G_1,\dots,G_k)$ for the set of all weakly join defined with respect to $G_1,\dots,G_k$.}
\end{defin}

\begin{theorem}\label{join}
Let $T_1$, $T_2$ and $T$  be  sets of transpositions in $S_n$ such that $G_T\in \mathcal{WJ}(G_{T_1}, G_{T_2})$. Suppose that $\Cay(S_n,T_1)$ is integral. Then $\Cay(S_n,T)$ is integral if and only if $\Cay(S_n,T_2)$ is integral.
\end{theorem}
\proof Let $T_3:=\{(i,i') \in S_n \mid \{i,i'\}\in E(G_T)\setminus (E(G_{T_1})\cup E(G_{T_2}))\}$. Consider the elements $a_i:=\sum_{t\in T_i} t$ ($i=1,2,3$)  in the group algebra $\mathbb{C}S_n$. It is easy to see that $\sum_{t\in T} t = a_1+a_2+a_3$ and in the group algebra $\mathbb{C}S_n$, the elements $a_1$, $a_2$ and $a_3$ are pairwise commuting.

Let $\rho_1,\dots, \rho_m$ be all inequivalent irreducible representations of $S_n$ and $d_i$ be the degree of $\rho_i$. Then, by Proposition \ref{at}, the spectrum of $\Cay(S_n,S)$ is the union of the spectra of $\rho_1(S),\dots,\rho_m(S)$, where $S\in\{T,T_1,T_2,T_3\}$.  Since $a_1$, $a_2$ and $a_3$ commute pairwise and $\rho_i$ is an algebra homomorphism, $\rho_i(a_j)$s commute pairwise. Hence there exists a basis $\mathcal{B}$ of the vector space $\mathbb{C}^{d_i}$ such that
the matrices $\rho_i(a_j)$, $j=1,2,3$ are diagonizable with respect to $\mathcal{B}$. We can assume that 
$\rho_i(a_j)=\diag(\lambda^i_{1j},\lambda^i_{2j},\dots,\lambda^i_{d_{i}j})$ with respect to the basis $\mathcal{B}$. On the other hand, since $\rho_i$ is an algebra homomorphism, we have $\rho_i(T)=\rho_i(a_1)+\rho_i(a_2)+\rho_i(a_3)$. This implies that the eigenvalues of  $\rho_i(T)$ are $\lambda^i_{k1}+\lambda^i_{k2}+\lambda^i_{k3}$, where $k=1,\dots,d_i$. By Proposition \ref{multipartite}, $\Cay(S_n, T_3)$ is integral. Thus $\Cay(S_n, T)$ is integral if and only if $\Cay(S_n,T_2)$ is integral.\eqd

\begin{remark}{\rm 
Let $T_1,\ldots,T_k$ be some sets of transpositions of $S_n$. Using Theorem \ref{join} and the induction on the number of graphs, one can easily see that if for every $i=1,\ldots,k$, the graph $\Cay(S_{n},T_i)$ is integral, then $\Cay(S_n,T)$ is integral, where $G_T\in \mathcal{WJ}(G_{T_1},\dots, G_{T_k})$.}
\end{remark}

\begin{remark}
{\rm Let $S=\{(1,2),(1,3),\ldots,(1,n)\}$ for some positive integer $n$. It is conjectured that
$\Cay(S_n,S)$ is integral in \cite{AV}. Using certain properties of the Jucys-Murphy elements, the validity of the latter conjecture was
confirmed  in \cite{CF}. Since the star graph $K_{1,n-1}\in \mathcal{WJ}(nK_1,nK_1)$, the latter conjecture is a direct consequence of Theorem \ref{join}.}
\end{remark}

A vertex $v$ of a graph $\Gamma$ is called a \textit{complete vertex} if it is adjacent to each vertex $u\neq v$ of $\Gamma$.
For a graph $\Gamma$ a \textit{cut vertex} of $\Gamma$ is a vertex $v$ of $\Gamma$ such that the number of connected components of $\Gamma-v$ is more than the number of connected components of $\Gamma$.
 As a direct consequence of the following result, we give a complete solution to Problem \ref{p5}.

\begin{corollary}\label{cut}
Let $T$ be a set of transpositions in $S_n$. Then 
\begin{itemize}
\item[(a)] If $G_T$ is a connected graph with a cut vertex $v$ then $\Cay(S_n,T)$ is integral if and only if $v$ is a complete vertex and $\Cay(S_{n},T')$ is integral, where $T'$ is a set of transpositions in $S_n$ such that $G_{T'}=(G_T-v)+v$.
\item[(b)] If $G_T$ is a tree, then $\Cay(S_n,T)$ is integral if and only if $G_T$ is the star graph $K_{1,n-1}$.
\end{itemize}
\end{corollary}
\proof 
$(a)$ Suppose that $\Cay(S_n,T)$ is integral. Then, by Theorem \ref{main} and \cite[Corollary 2.1]{K}, $v$ is a complete vertex.
We know that $T''$ is the set of transpositions in $S_n$ such that $G_{T''}=nK_1$ then $T''=\varnothing$ and $\Cay(S_n,T'')$ is integral.
Thus as $G_T\in\mathcal{WJ}(G_{T'},nK_1)$, by Theorem \ref{join}, $\Cay(S_{n},T')$ is integral. This proves one direction. Conversely, suppose that $\Cay(S_n,T')$ is integral and $v$ is a complete vertex. Then again as $G_T\in\mathcal{WJ}(G_{T'},nK_1)$, Theorem \ref{join} implies that $\Cay(S_n,T)$ is integral.

$(b)$ It is clear that if $G_T=K_1$ or $K_2$ there is nothing to prove. So we may assume that $G_T$ has at least three vertices. Suppose that $\Cay(S_n,T)$ is integral. Then as $G_T$ is a tree, by part $(a)$, $G_T$ has a complete vertex and so it is $K_{1,n-1}$. The converse direction is clear by Proposition \ref{multipartite}.
\eqd

\begin{defin} \label{def1} {\rm   a) Let $\Gamma_1,\dots,\Gamma_k$ be graphs with the same vertex set and for each $i$, $\Gamma_i$ be a disjoint union of complete multipartite graphs. Then a graph $\Gamma$ is called \textit{generalized complete multipartite graph of  type $1$} if $\Gamma\in\mathcal{WJ}(\Gamma_1,\dots,\Gamma_k)$.\\
b) For $n\geq 2$, we say that the graph $\Gamma$ is a \textit{generalized complete multipartite graph of  type $n$}, if $\Gamma\in\mathcal{WJ}(\Gamma_1,\dots,\Gamma_k)$, where $\Gamma_1,\dots,\Gamma_k$ are graphs with the same vertex set and for every $i=1,2,\dots k$, the graph $\Gamma_i$ is a disjoint union of generalized complete multipartite graphs of  type at most $n-1$.}
\end{defin}

\begin{proposition}\label{typek}
Let $T$ be a set of transpositions in $S_n$ such that $G_T$ is a generalized complete multipartite graph of type $k$, for some positive integer $k$. Then $\Cay(S_n,T)$ is integral.
\end{proposition}
\proof We prove by induction on $k$. The case $k=1$ is clear by Proposition \ref{multipartite} and Theorem \ref{join}. Suppose that the result is true for every generalized complete multipartite graph of type at most $k-1$. Now the result follows from the induction hypothesis, \cite[Proposition 5.2]{C} and Theorem \ref{join}.
\eqd

By the above result in the sequel, we suggest the following conjecture and try to show that this conjecture is true for some classes of Cayley graphs of symmetric groups generated by transpositions.

\begin{conj}\label{conj} Let $T$ be a set of transpositions in $S_n$. Then $\Cay(S_n, T)$ is integral, if and only if for some positive integer $k$,  the graph $G_T$ is a generalized complete multipartite graph of  type $k$.
\end{conj}

The following  results confirm Conjecture \ref{conj} for cycles.

\begin{proposition}\label{prop1}
Let $T$ be a set of transposition in $S_m$ such that $G_T\cong C_m$. Then $\Cay(S_m,T)$ is integral if and only if $m=3$ or $4$.
\end{proposition}
\proof It is well-known that the only integral cycles are $C_3$, $C_4$ and $C_6$, see for example \cite[p. 43]{BCRSS}. Since cycles are regular graphs, the above cycles are the only Laplacian integral cycles. Thus, one direction is a direct consequence of 
Theorem \ref{main} and Remark \ref{rmk}. The converse direction is clear by Proposition \ref{multipartite}.\eqd 

\begin{proposition} \label{prop2} Let $T$ be a set of transpositions in $S_n$ such that $G_T$ is a connected cubic graph. Then $\Cay(S_n,T)$ is integral if and only if $G_T$ is isomorphic to one of the graphs $K_4$ or $K_{3,3}$.
\end{proposition}
\proof Since $G_T$ is a connected cubic integral graph, $G_T$ is one of the thirteen graphs given in \cite{schwenk}.  Suppose that $\Cay(S_n,T)$ is integral.  Then $G_T$ is isomorphic to $K_4$ or $K_{3,3}$, by GAP \cite{GAP2022}, Proposition \ref{at} and Theorem \ref{main}. The converse direction is clear by Proposition \ref{multipartite}.\eqd

A connected graph with $n$ vertices and $m$ edges is called a \textit{$k$-cyclic graph} if $k=m-n+1$. In the following result, we give a characterization of integral Cayley graphs $\Cay(S_n,T)$, where $G_T$ is a $k$-cyclic graph and $k=1,2,3$ which confirms Conjecture \ref{conj}.

\begin{proposition}\label{prop3}
Let $T$ be a set of transpositions in $S_n$ such that $G_T$ is a $k$-cyclic graph.
\begin{itemize}
\item[(a)] If $k=1$ then $\Cay(S_n,T)$ is integral if and only if  $G_T$ is one of the graphs $K_3$, $K_{2,2}$ or $(K_2\cup (n-3)K_1)*K_1$.
\item[(b)] If $k=2$ then $\Cay(S_n,T)$ is integral if and only if $G_T$ is one of the graphs $K_{2,3}$, $K_2*(2K_1)$, $(2K_2\cup (n-5)K_1)*K_1$ or $(K_{1,2}\cup (n-4)K_1)*K_1$.
\item[(c)] If $k=3$ then $\Cay(S_n,T)$ is integral if and only if $G_T$ is one of the graphs $(3K_2)*K_1$, $(K_{1,2}\cup K_2)*K_1$, $K_{2,4}$, $K_{1,3}*K_1$, $K_4$, $K_{2,3}$, $(3K_2\cup (n-7)K_1)*K_1$, $(K_{1,2}\cup K_2\cup (n-6)K_1)*K_1$, $(K_3\cup (n-4)K_1)*K_1$ or $(K_{1,3}\cup (n-5)K_1)*K_1$.
\end{itemize}
\end{proposition}
\proof
 Suppose that $\Cay(S_n,T)$ is integral. Then, by Theorem \ref{main}, $G_T$ is Laplacian integral. If $k=1$, then \cite[Theorem 3.2]{LL} and Remark \ref{rmk} imply that $G_T$ is isomorphic to one of the graphs  $K_3$, $K_{2,2}$ or $(K_2\cup (n-3)K_1)*K_1$.

If $k=2$ then by \cite[Theorem 3.3]{LL}, $G_T$ is one of the graphs $(2K_2\cup (n-5)K_1)*K_1$, $(K_{1,2}\cup (n-4)K_1)*K_1$, $K_{2,3}$, $K_2*(2K_1)$ or $K_{1,2}\square K_2$, where $\square$ is the direct product operation of graphs. In the last case, one can easily see that $\Cay(S_6,T)$ is not integral, by GAP \cite{GAP2022} and Proposition \ref{at}.

If $k=3$ then by \cite[Theorem 4.1 and Theorem 5.1]{HHW}, $G_T$ is one of the thirteen graphs $(3K_2)*K_1$, $(K_{1,2}\cup K_2)*K_1$, $K_{2,4}$, $K_{1,3}*K_1$, $K_4$, $K_{2,3}$, $(3K_2\cup (n-7)K_1)*K_1$, $(K_{1,2}\cup K_2\cup (n-6)K_1)*K_1$, $(K_3\cup (n-4)K_1)*K_1$, $(K_{1,3}\cup (n-5)K_1)*K_1$, or $\Gamma_1$, $\Gamma_2$ and $\Gamma_3$ shown in Fig. 1. We see, by GAP \cite{GAP2022} and Proposition \ref{at}, that in the last three cases $\Cay(S_n,T)$ is not integral. The converse direction is clear by Proposition \ref{typek}. This completes the proof.\eqd

\begin{figure}\label{tricyclic}
\centering{\includegraphics[scale=.7]{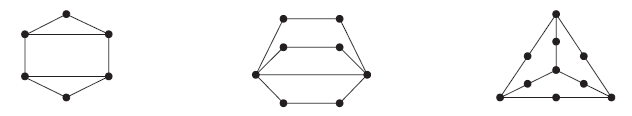}}
\caption{\textit{The Laplacian integral tricyclic graphs $\Gamma_1$, $\Gamma_2$ and $\Gamma_3$.} }
\end{figure}

\section*{Acknowledgments}
The work of Mahdi Ebrahimi is supported in part
by a grant  from School of Mathematics, Institute for Research in Fundamental Sciences (IPM).

\end{document}